\documentclass{svproc}
\usepackage{url}

\usepackage{amsfonts}

\usepackage{latexsym}
\usepackage{amssymb}
\usepackage{amsmath}
\usepackage[mathscr]{eucal}
\usepackage{graphicx}
\usepackage{hyperref}
\usepackage{caption}
\usepackage{subcaption}

\renewcommand{\qed}{\hfill{\ \ \rule{2mm}{2mm}} \vspace{0.2in}}
\newcommand{\ind}{1\hspace{-2.3mm}{1}}

\begin{document}
\mainmatter              
\title{Edge Constrained Eulerian Extensions}
\titlerunning{Edge Constrained Eulerian Extensions}  
%
\author{Ghurumuruhan Ganesan \thanks{Corresponding Author}}
\authorrunning{G. Ganesan}

%
%
\institute{IISER Bhopal,\\\email{gganesan82@gmail.com}}

\maketitle              

\begin{abstract}
In this paper we study Eulerian extensions with edge constraints and use the   probabilistic method to establish sufficient conditions for a given connected graph to be a subgraph of a Eulerian graph containing~\(m\) edges, for a given number~\(m.\)

\keywords{Eulerian Extensions, Edge Constraint, Probabilistic Method}
\end{abstract}
\renewcommand{\theequation}{\thesection.\arabic{equation}}
\setcounter{equation}{0}
\section{Introduction} \label{intro}

In the Eulerian extension problem, a given graph is to be converted into an Eulerian graph by addition of as few edges as possible and such problems have applications in routing and scheduling (Dorn et al. (2013)).   Boesch et al. (1977) studied conditions under which a graph~\(G\) can be extended to an Eulerian graph and later Lesniak and Oellermann (1986) presented a detailed survey on subgraphs and supergraphs of Eulerian graphs and multigraphs. For applications of Eulerian extensions to scheduling and parametric aspects, we refer to H\"ohn et al. (2012) and  Fomin and Golovach (2012), respectively.



In this paper, we construct Eulerian extension of graphs with a predetermined number of edges. Specifically, given a graph~\(G\) with maximum degree~\(\Delta\) and~\(b\) number of edges and given an integer~\(m > b,\) we use the asymmetric probabilistic method to derive sufficient conditions for the existence of an Eulerian extension of~\(G\) with~\(m\) edges.

The paper is organized as follows: In Section~\ref{sec_eul_ext}, we state and prove our main result regarding Eulerian extensions with edge constraints.




\renewcommand{\theequation}{\thesection.\arabic{equation}}
\setcounter{equation}{0}
\section{Edge Constrained Eulerian Extensions}\label{sec_eul_ext}
Let~\(G = (V,E)\) be a graph with vertex set~\(V\) and edge set~\(E.\) The vertices~\(u\) and~\(v\) are said to be adjacent in~\(G\) if the edge~\((u,v)\) with endvertices~\(u\) and~\(v\) is present in~\(E.\) We define~\(d_G(v)\) to be the degree of vertex~\(v,\) i.e., the number of vertices adjacent to~\(v\)
in~\(G.\)

A sequence of vertices~\({\cal W} := (u_1,u_2,\ldots,u_t)\) is said to be a \emph{walk}  if~\(u_i\) is adjacent to~\(u_{i+1}\) for each~\(1 \leq i \leq t-1.\) If in addition the vertex~\(u_t\) is also adjacent to~\(u_1,\) then~\({\cal W}\) is said to be a \emph{circuit}. We say that~\({\cal W}\) is an Eulerian circuit if each edge of the graph~\(G\) occurs exactly once in~\({\cal W}.\) The graph~\(G\) is said to be an \emph{Eulerian} graph if~\(G\) contains an Eulerian circuit. 



Let~\(G\) be any graph. We say that a graph~\(H\) is an \emph{Eulerian extension} of~\(G\) if~\(G\) and~\(H\) share the same vertex set,~\(G\) is a subgraph of~\(H\) and~\(H\) is Eulerian.
\begin{definition}\label{def_one}
For an integer~\(m \geq 1,\) we say that a graph~\(G\) is~\(m\)-\emph{Eulerian extendable} if there exists an Eulerian extension~\(H\) of~\(G\) containing exactly~\(m\) edges.
\end{definition}

We have the following result regarding~\(m\)-Eulerian extendability. Throughout, constants do not depend on~\(n.\)
\begin{theorem}\label{thm_comp}
For every pair of constants~\(0 < \alpha,\beta < 1\) satisfying~\(\beta + 40\alpha^2 < \frac{1}{2}\) strictly, there exists a constant~\(N = N(\alpha,\beta) \geq 1\) such that the following holds for all~\(n \geq N\): Let~\(m\) be any integer satisfying
\begin{equation}\label{m_range2}
2n \leq m  \leq \alpha \cdot n^{\frac{3}{2}}
\end{equation}
and let~\(G \subset K_n \) be any connected graph containing~\(n\) vertices,~\(b\) edges and a maximum vertex degree~\(\Delta.\) If
\begin{equation}\label{m_cond2}
\Delta \leq \beta \cdot n \text{ and }b \leq m-n,
\end{equation}
then~\(G\) is~\(m\)-Eulerian extendable.
\end{theorem}
To see the necessity of the bound~\(b \leq m-n,\) we use the fact that a graph~\(H\) is Eulerian if and only if~\(H\) is connected and each vertex of~\(H\) has even degree (Theorem 1.2.26, pp. 27, West (2001)). Therefore to obtain an Eulerian extension of~\(G,\) we only need to convert all odd degree vertices into even degree vertices.

Suppose that~\(n\) is even and all the vertices in~\(G\) have an odd degree. Because the degree of each vertex is at most~\(\frac{n}{2}-1\) (see~(\ref{m_cond2})) the sum of neighbourhood sizes of~\(2i-1\) and~\(2i\) is at most~\(n-2.\) Therefore for each~\(1 \leq i \leq \frac{n}{2},\) there exists a vertex~\(w_i\) neither adjacent to~\(2i-1\) nor adjacent to~\(2i\) in  the graph~\(G.\) Adding the~\(n\) edges~\(\{(w_i,2i-1),(w_i,2i)\}_{1 \leq i \leq \frac{n}{2}}\)  gives us an Eulerian extension of~\(G.\)

In our proof of Theorem~\ref{thm_comp} below, we use the asymmetric probabilistic method for higher values of~\(m\) to obtain walks of predetermined lengths between pairs of odd degree vertices and thereby construct the desired extension,





\subsection*{Proof of Theorem~\ref{thm_comp}}
As before, we use the fact that a graph~\(H\) is Eulerian if and only if~\(H\) is connected and each vertex of~\(H\) has even degree. We assume that the vertex set of~\(G\) is~\(V := \{0,1,2,\ldots,n-1\}\) and also let~\({\cal T}\) be the set of all odd degree vertices in~\(G\) so that the number of odd degree vertices~\(\#{\cal T},\) is even.  If there are vertices~\(u,v \in {\cal T}\) that are not adjacent to each other in~\(G,\) then we mark the edge~\((u,v)\) and also the endvertices~\(u\) and~\(v.\) We then pick two new non-adjacent vertices~\(x\) and~\(y\) in~\({\cal T} \setminus \{u,v\}\) and repeat the procedure. We continue this process until we reach one of the following two scenarios: Either the number of marked edges is~\(m-b\) in which case, we simply add the marked edges to~\(G\) and get the desired Eulerian extension~\(H.\) Or, we are left with a set of marked edges of cardinality, say~\(l\) and a clique~\({\cal C} := \{u_1,\ldots,u_{2z}\} \subset {\cal T}\) containing~\(2z\) unmarked vertices.


Let~\(G_0\) be the graph obtained by adding all the~\(l \leq \frac{n}{2}\) marked edges to~\(G.\) If~\(\Delta_0\) and~\(b_0\) denote the maximum vertex degree and the number of edges in~\(G_0,\) respectively, then~
\begin{equation}\label{del_not}
\Delta_0 \leq \Delta+1 \text{ and }b_0 =b+l \leq b + n.
\end{equation}
We now pair the vertices in~\({\cal C}\) as~\(\{u_{2i-1},u_{2i}\}_{1 \leq i \leq z}\) assuming that~\(z \geq 1\) (If not, we simply remove a marked edge~\(e\) from~\(G_0\) and label the endvertices of~\(e\) as~\(u_1\) and~\(u_2\)). We use the probabilistic method to obtain~\(z\) edge-disjoint walks~\(\{{\cal W}_i\}_{1 \leq i \leq z}\) containing no edge of~\(G_0\) such that each walk~\({\cal W}_i\) has~\(w\) edges and~\(u_{2i-1}\) and~\(u_{2i}\) as endvertices, where~\(w\) satisfies
\begin{equation}\label{w_def}
b_0 + z\cdot w = m.
\end{equation}
Adding the walks~\(\{{\cal W}_i\}_{1 \leq i \leq z}\) to~\(G_0\) would then give us the desired~\(m\)-Eulerian extension.

In~(\ref{w_def}) we have assumed for simplicity that~\(w = \frac{m-b_0}{z}\) is an integer. If not, we write~\(m-b_0 = z \cdot w + r \) where~\(0 \leq r \leq w-1\) and construct the~\(z-1\) walks~\({\cal W}_i, 1 \leq i \leq z-1\) each of length~\(w\) edges and the last walk~\({\cal W}_z\) of length~\(w+r \leq 2w.\) Again adding these walks to~\(G_0\) would give us the desired Eulerian extension with~\(m\) edges.
For future use we remark that the length~\(w\) of each walk added in the above process is bounded above by
\begin{equation}\label{w_up}
w = \frac{m-b_0}{z} \leq m \leq \alpha \cdot n^{\frac{3}{2}}.
\end{equation}



We begin with the pair of vertices~\(u_1\) and~\(u_2.\) Let~\(\{X_i\}_{1 \leq i \leq w}\) be independent and identically distributed (i.i.d.) random variables uniformly distributed in the set~\(\{0,1,\ldots,n-1\}.\) Letting~\({\cal S} := (u_1,X_1,\ldots,X_{w},u_2),\) we would like to convert the sequence~\({\cal S}\) into a walk~\({\cal W}_1\) with endvertices~\(u_1\) and~\(u_2\) and containing no edge of~\(G_0.\) The construction of~\({\cal W}_1\) is split into two parts: In the first part, we collect the preliminary relevant properties of~\({\cal S}\) and in the second part, we obtain the walk~\({\cal W}_1.\) \\\\
\emph{\underline{Preliminary definitions and estimates}}: An entry in~\({\cal S}\) is defined to be a vertex and we define~\((u_1,X_1),(X_w,u_2)\) and~\(\{(X_i,X_{i+1})\}_{1 \leq i \leq w-1}\) to be the edges of~\({\cal S}.\) The neighbour set of a vertex~\(v\) in~\({\cal S}\) the set of vertices~\(u\) such that either~\((v,u)\) or~\((u,v)\) appears as an edge of~\({\cal S}.\) The neighbour set of~\(v\) in the multigraph~\(G_0 \cup {\cal S}\) is the union of the neighbour set of~\(v\) in the graph~\(G_0\) and the neighbour set of~\(v\) in~\({\cal S}.\) The degree of a vertex~\(v\) in~\(G_0 \cup {\cal S}\) is defined to be the sum of the degree of~\(v\) in~\(G_0\) and the degree of~\(v\) in~\({\cal S}.\)

The three main ingredients used in the construction of the walk~\({\cal W}_1\) are:\\
\((1)\) The degree of a vertex in the multigraph~\(G_0 \cup {\cal S},\)\\
\((2)\) the number of ``bad" vertices in~\(G_0 \cup {\cal S}\) and\\
\((3)\) the number of ``bad" edges in~\({\cal S}.\)\\
Below, we define and estimate each of the three quantities in that order.

We first estimate the degree of each vertex in the multigraph~\(G_0 \cup {\cal S}.\) For any~\(0 \leq v \leq n-1\) and any~\(1 \leq i \leq w,\) let~\(I_i = \ind(X_i = v)\) be the indicator function of the event that~\(X_i =v.\) We have~\(\mathbb{P}(I_i=1) = \frac{1}{n}\) and so if~\(D_v = \sum_{i=1}^{w} I_i\) denotes the number of times the entry~\(v\) appears in the sequence~\((X_1,\ldots,X_{w}),\) then~\(\mathbb{E}D_v = \frac{w}{n}\) and so by the standard deviation estimate~(\ref{conc_est_f}) in Appendix, we have
\begin{equation}\label{dv_est}
\mathbb{P}\left(D_v \geq \frac{2w}{n} \right)\leq 2\exp\left(-\frac{w}{16n} \right).
\end{equation}
If~\(\frac{w}{n} \geq 100 \log{n}\) then we get from~(\ref{dv_est}) that~\(\mathbb{P}(D_v \geq \frac{2w}{n}) \leq \frac{1}{n^2}.\)
Else, we use Chernoff bound directly to get that
\begin{equation}\label{dv_est_p}
\mathbb{P}\left(D_v \geq 100\log{n} \right)\leq \frac{1}{n^2}.
\end{equation}
Therefore setting~\(a_n := \max\left(\frac{2w}{n},100\log{n}\right),\) we get that
\begin{equation}\label{dv_est_pp}
\mathbb{P}\left(D_v \geq a_n \right)\leq \frac{1}{n^2}.
\end{equation}

If the event
\begin{equation}\label{e_deg_def}
E_{deg} := \bigcap_{0 \leq v \leq n-1} \left\{D_v \leq a_n\right\}
\end{equation}
occurs, then in~\(G_0 \cup {\cal S}\) each vertex has degree at most~\(\Delta_0+1+a_n,\) with the extra term~\(1\) to account for the fact that vertices~\(X_1\) and~\(X_w\) are also adjacent to~\(u_1\) and~\(u_2,\) respectively. By the union bound and~(\ref{dv_est}) we therefore have
\begin{equation}\label{first_cond}
\mathbb{P}(E_{deg}) \geq 1-\frac{1}{n}.
\end{equation}

The next step is to estimate the number of ``bad" vertices in~\(G_0 \cup {\cal S}.\)  Let~\(X_0 := u_1, X_{w+1} := u_2\) and for~\(0 \leq i \leq w-1,\) say that vertex~\(X_i\) is \emph{bad} if~\(X_i = X_{i+1}\) or~\(X_i = X_{i+2}.\) For simplicity define~\(X_w\) to be bad always. If~\(J_i\) is the indicator function of the event that vertex~\(X_i\) is bad, then for~\(0 \leq i \leq w-1,\) we have that
\begin{equation}\label{z_est}
\frac{1}{n} \leq \mathbb{P}(J_i=1) \leq \frac{2}{n}.
\end{equation}
The term~\(N_{v,bad} := \sum_{i=0}^{w-1}J_i+1\) denotes the total number of bad vertices in the sequence~\({\cal S}.\) To estimate~\(N_{v,bad}\) we split~\(N_{v,bad}-1   = J(A) + J(B) + J(C),\) where~\[J(A) = J_1 + J_4 + \ldots, J(B) = J_2 + J_5 + \ldots \text{ and } J(C) = J_3 + J_6 + \ldots\] so that each~\(J(u), u \in \{A,B,C\}\) is a sum of i.i.d.\ random variables.

The term~\(J(A)\) contains at least~\(\frac{w}{3}-1\) and at most~\(\frac{w}{3}\) random variables. As in the proof of~(\ref{dv_est_pp}), we use~(\ref{z_est}) and the standard deviation estimate~(\ref{conc_est_f}) in Appendix to obtain that
\[\mathbb{P}\left(J(A) \geq \frac{a_n}{3} \right) \leq \frac{1}{n^2}\] for all~\(n\) large. A similar estimate holds for~\(J(B)\) and~\(J(C)\) and so combining these estimates and using the union bound, we get that
\begin{equation}\label{second_cond}
\mathbb{P}\left(E_{v,bad}\right) \geq 1-\frac{3}{n^2}
\end{equation}
where~\(E_{v,bad} := \{N_{v,bad}\leq a_n +1\}\) denotes the event that the number of bad vertices in~\({\cal S}\) is at most~\(a_n+1.\)

The final estimate involves counting the number of bad edges in the sequence~\({\cal S}.\) For~\(0 \leq i \leq w\) say that~\((X_{i},X_{i+1})\) is a \emph{bad edge} if one of the following two conditions hold:\\
\((d1)\) Either~\(\{X_{i},X_{i+1}\}\) is an edge of~\(G_0\) or\\
\((d2)\) There exists~\(i+2 \leq j \leq w\) such that~\(\{X_{i},X_{i+1}\} = \{X_{j},X_{j+1}\}.\)\\
To estimate the probability of occurrence of~\((d1),\) let~\(e\) be an edge of~\(G_0\) with endvertices~\(u\) and~\(v.\) We have that
\[\mathbb{P}\left(\{X_{i},X_{i+1}\} = \{u,v\}\right) \leq \frac{2}{n^2}.\] Similarly for any~\(i+2 \leq j \leq w,\) the possibility~\((d2)\) also occurs with probability at most~\(\frac{2}{n^2}.\) Therefore if~\(L_i\) is the indicator function of the event that~\((X_i,X_{i+1})\) is a bad edge, we have that
\begin{equation}\label{ji_est}
\mathbb{P}(L_{i}=1) \leq \sum_{l=1}^{b_0} \frac{2}{n^2}  + \sum_{j=i+2}^{w}\frac{2}{n^2} \leq \frac{2(b_0+w)}{n^2}.
\end{equation}

If~\(N_{e,bad} := \sum_{i=0}^{w} L_i\) denotes the total number of bad edges in~\({\cal S},\) then from~(\ref{ji_est}) and the fact that~\(L_0 \leq 1\) we have
\begin{equation}\label{cn_def}
\mathbb{E}N_{e,bad} \leq 1+\frac{2(b_0+w)w}{n^2} =:c_n.
\end{equation}
Letting~\(E_{e,bad} := \{N_{e,bad} \leq K \cdot c_n\}\) denote the event that the number of bad edges in~\({\cal S}\) is at most~\(K \cdot c_n,\)
for some large integer constant~\(K \geq 1\) to be determined later, we get from Markov inequality that
\begin{equation}\label{third_cond}
\mathbb{P}\left(E_{e,bad}\right) \geq 1-\frac{1}{K},
\end{equation}

If~\(E_{valid}\) denotes the event that the first and last edges~\((X_0,X_1)\) and~\((X_{w},X_{w+1})\) are valid edges not in~\(G_0,\) then using the fact that the degree of any vertex in~\(G_0\) is at most~\(\frac{n}{2}\) (see~(\ref{del_not}) and~(\ref{m_cond2}) in the statement of the Theorem) we get that
\begin{equation}\label{valid_cond}
\mathbb{P}(E_{valid}) \geq \left(\frac{1}{2}-\frac{1}{n}\right)^2 . 
\end{equation}
Defining the joint event
\[E_{joint}:= E_{valid} \cap E_{deg} \cap E_{v,bad} \cap E_{e,bad}\]
and using
\begin{eqnarray}
\mathbb{P}\left(A \bigcap \bigcap_{i=1}^{l} B_i\right) &\geq& \mathbb{P}(A) - \mathbb{P}\left(\bigcup_{i=1}^{l} B_i^c\right) \nonumber\\
&\geq& \mathbb{P}(A) - \sum_{i=1}^{l}\mathbb{P}\left(B_i^c\right). \label{asym_prob}
\end{eqnarray}
with~\(A = E_{valid}\) we get from~(\ref{first_cond}),~(\ref{second_cond}),~(\ref{third_cond}) and~(\ref{valid_cond}) that
\begin{equation}\label{a_joint_est}
\mathbb{P}(E_{joint}) \geq \left(\frac{1}{2}-\frac{1}{n}\right)^2-\frac{1}{K}- \frac{1}{n} -\frac{3}{n^2} \geq \frac{1}{21}
\end{equation}
for all~\(n\) large, provided the constant~\(K = 5,\) which we fix henceforth. This completes the preliminary estimates used in the construction of the walk~\({\cal W}_1.\)\\\\
\emph{\underline{Construction of the walk~\({\cal W}_1\)}}: Assuming that the event~\(E_{joint}\) occurs, we now convert~\({\cal S}_0 := {\cal S}\) into a walk~\({\cal W}_1.\) We begin by ``correcting" all bad vertices. Let~\(X_{i_1},X_{i_2},\ldots,X_{i_t}, i_1 < i_2 < \ldots <i_t\) be the set of all bad vertices. Thus for example either~\(X_{i_1} = X_{i_1+1}\) or~\(X_{i_1} = X_{i_1+2}.\) Because the event~\(E_{deg}\) occurs, we get from the discussion following~(\ref{e_deg_def}) that the degree of each vertex in~\(G_0 \cup {\cal S}_0\) is at most~\(\Delta_0+a_n +1.\) From~(\ref{del_not}) and the first condition in~(\ref{m_cond2}) we get that
\begin{equation}\label{del_not_est}
\Delta_0 \leq \Delta+1 \leq \frac{n}{3}+1
\end{equation}
and from the definition of~\(a_n\) prior to~(\ref{dv_est_pp}) and the upper bound~\(w \leq n^{\frac{3}{2}}\) in~(\ref{w_up}), we get that
\begin{equation}\label{an_est}
a_n = \max\left(100\log{n}, \frac{2w}{n}\right)\leq 100\log{n} + \frac{2w}{n} \leq 100 \log{n} + 2\sqrt{n} \leq 3\sqrt{n}
\end{equation}
for all~\(n\) large. Consequently, using~\(\beta< \frac{1}{2}\) strictly (see statement of Theorem~\ref{thm_comp}) ,
\begin{equation}
\Delta_0+a_n +1 \leq \beta \cdot n  + 1 + 3\sqrt{n}  \leq \frac{n}{2}-5\label{del_a}
\end{equation}
for all~\(n\) large. From~(\ref{del_a}), we therefore get that there exists a vertex~\(v_1\) that is \emph{not} a neighbour of~\(X_{i_1}\) in~\(G_0 \cup {\cal S}.\) Similarly, the total number of neighbours of~\(v_1\) and~\(X_{i_1+3}\) in~\(G_0 \cup {\cal S}\) is at most
\begin{equation}\label{deg_tot}
2\Delta_0  + 2a_n+2 \leq 2\beta \cdot n + 2 + 6\sqrt{n}  < n-10
\end{equation}
for all~\(n\) large and so there exists a vertex~\(v_2 \neq X_{i_1}\) that is not a neighbour of~\(v_1\) and also not a neighbour of~\(X_{i_1+3}\) in~\(G_0 \cup {\cal S}.\)

We now set~\(X^{(1)}_{i_1+1} = v_1\) and~\(X^{(1)}_{i_1+2} = v_2\)  and~\(X^{(1)}_{j} = X_j\) for~\(j \neq i_1+1,i_1+2\) and call the resulting sequence as~\({\cal S}_1 := (X^{(1)}_1,\ldots,X^{(1)}_{w}).\) By construction the degree of each vertex in the multigraph~\(G_0 \cup {\cal S}_1\) is at most~\(\Delta+a_n+1+2\) and there are at most~\(t-1\) bad vertices in~\({\cal S}_1.\) We now pick the bad vertex with the least index in~\({\cal S}_1\) and repeat the above procedure with~\({\cal S}_1\) to get a sequence~\({\cal S}_2\) containing at most~\(t-2\) bad vertices.

After~\(k \leq t\) iterations of the above procedure, the degree of each vertex in the multigraph~\(G_0 \cup {\cal S}_k\) would be at most
\begin{equation}\label{st_def}
\Delta_0+a_n + 1+2k \leq \Delta_0+a_n +1+2t \leq \Delta_0+3a_n + 3
\end{equation}
because the event~\(E_{joint} \subseteq E_{v,bad}\) occurs and so~\(t \leq a_n+1.\) Again using~(\ref{del_not_est}) and~(\ref{an_est}) and arguing as in~(\ref{deg_tot}), we get that the sum of the degrees of any two vertices in~\(G_0 \cup {\cal S}_t\) is at most~\(n-10\) for all~\(n\) large. Thus the above procedure indeed proceeds for~\(t\) iterations and by construction, the sequence~\({\cal S}_t\) obtained at the end, has no bad vertices.

We now perform an analogous procedure for correcting all bad edges in~\({\cal S}_t.\) For example if~\((X_{l},X_{l+1})\) is a bad edge in~\({\cal S}_t,\) then following an analogous argument as before we pick a vertex~\(Y_{l+1}\) that is neither adjacent to~\(X_l\) nor adjacent to~\(X_{l+2}\) in the sequence~\({\cal S}_t.\) We replace~\(X_{l+1}\) with~\(Y_{l+1}\) to get a new sequence~\({\cal S}_{t+1}.\) In the union~\(G_0 \cup {\cal S}_{t+1}\) the degree of each vertex is at most~\(\Delta_0+3a_n+3+2\) (see~(\ref{st_def})) and the number of bad edges is at most~\(r-1.\) At the end of~\(r \leq K \cdot c_n\) iterations, we obtain a multigraph~\(G_0 \cup {\cal S}_{t+r},\) where the degree of each vertex is at most
\[\Delta_0 + 3a_n + 3+2r \leq \Delta_0 + 3a_n+3+2Kc_n,\] since the event~\(E_{e,bad}\) occurs  and therefore~\(N_{e,bad} \leq K \cdot c_n\) (see discussion preceding~(\ref{third_cond})). Substituting the expression for~\(c_n\) from~(\ref{cn_def}) and using the second estimate for~\(a_n \) in~(\ref{an_est}), we get that~\(\Delta_0 + 3a_n + 3+2Kc_n\) is at most
\begin{eqnarray}
&&\Delta_0 + 300\log{n}+ \frac{6w}{n} + 3+2K+ \frac{2K(b_0+w)w}{n^2} \nonumber\\
&&\;\;\;\;\leq\;\;\;\Delta_0 + 301\log{n} + \frac{6w}{n} + \frac{2K(b_0+w)w}{n^2} \label{con_est}
\end{eqnarray}
for all~\(n\) large. Recalling that~\(u_1\) and~\(u_2\) are the endvertices of the starting sequence~\({\cal S}_0,\) we get that the final sequence~\({\cal S}_{t+r}\) contains no bad edge and is therefore the desired walk~\({\cal W}_1\) with endvertices~\(u_1\) and~\(u_2.\) This completes the construction of the walk~\({\cal W}_1.\)\\\\
\emph{\underline{Rest of the walks}}: We now repeat the above procedure to construct the rest of the walks. We set~\(G_1 := G_0 \cup {\cal W}_1\) and argue as above to obtain a walk~\({\cal W}_2\) with~\(w\) edges present in~\(\overline{G}_1\) and containing~\(u_3\) and~\(u_4\) as endvertices. Adding the walk~\({\cal W}_2\) to~\(G_1\) we get a new graph~\(G_2.\) In effect, to the graph~\(G_1\) containing~\(b_0+w\) edges, we have added~\(w\) edges and by an argument analogous to~(\ref{con_est}), we have increased the degree of a vertex by at most~\[301\log{n} + \frac{6w}{n} + \frac{2K(b_0+2w)w}{n^2},\] in obtaining the graph~\(G_2.\) We recall that (see first paragraph of the proof) there are~\(z\) such walks to be created of which~\(z-1\) have length~\(w\) and the final walk has length at most~\(2w.\) Therefore after~\(z\) iterations, we get a graph~\(G_z\) with~\(m\) edges and whose maximum vertex degree~\(\Delta_z\) is at most
\begin{eqnarray}\label{del_z_est}
\Delta_z &\leq&\Delta_0 + \left(301\log{n} + \frac{6w}{n}\right) \cdot (z-1) + 301\log{n}+\frac{12w}{n} \nonumber\\
&&\;\;\;\;+\;\;2K\sum_{k=1}^{z-1}\frac{(b_0+k\cdot w)w }{n^2} + 2K\frac{(b_0+(z-1)\cdot w)2w }{n^2}.
\end{eqnarray}
By construction~\(G_z\) is an Eulerian graph. 

To verify the obtainability of~\(G_z,\) we estimate~\(\Delta_z\) as follows. The term~\(z\) is no more than the size of a maximum clique in the original graph~\(G\) (see discussion prior to~(\ref{del_not})) and since there are~\(m \leq n^{\frac{3}{2}}\) edges in~\(G,\) the maximum size of a clique in~\(G\) is at most~\(n^{\frac{3}{4}}.\) Therefore
\begin{equation}\label{z_bound}
z \leq n^{\frac{3}{4}}.
\end{equation}
Also using~(\ref{w_def}) and~(\ref{m_cond2}), we get that~\(zw \leq m \leq \alpha \cdot n^{\frac{3}{2}}\) and so
\begin{equation}\label{wz_bound}
\frac{wz}{n} \leq \alpha \cdot \sqrt{n} \leq \sqrt{n}
\end{equation}
Finally from~(\ref{del_not}) we have that~\(b_0 \leq b+n\) and so the second line in~(\ref{del_z_est}) is at most
\begin{eqnarray}
\sum_{k=1}^{z} \frac{(b+n+k\cdot w)w }{n^2} &\leq& \frac{z(b+n+zw)2w}{n^2}  \nonumber\\
&\leq& \frac{2m(n+m)}{n^2} \nonumber\\
&\leq& 2\sqrt{n} + \frac{2m^2}{n^2} \nonumber\\
&\leq& \sqrt{n} + 2\alpha^2 \cdot n \label{term_3}
\end{eqnarray}
where the second inequality in~(\ref{term_3}) follows from the estimate~\(b+zw \leq b_0+zw=m\) (see~(\ref{w_def})), the third and fourth estimates in~(\ref{term_3}) follow from the bound~\(m \leq \alpha \cdot n^{\frac{3}{2}}\) (see~(\ref{m_range2})).

Plugging~(\ref{term_3}),~(\ref{wz_bound}) and~(\ref{z_bound}) into~(\ref{del_z_est}) we get that
\begin{eqnarray}
\Delta_z &\leq& \Delta_0 + 301n^{\frac{3}{4}} \cdot \log{n} + \sqrt{n}\left(\frac{12}{10} + 2K\right) + 4K\alpha^2 \cdot n \nonumber\\
&\leq& (\beta + 4K\alpha^2) \cdot n + 1 + 301n^{\frac{3}{4}} \cdot \log{n} + \sqrt{n}\left(\frac{12}{10} + 2K\right)\label{tata}
\end{eqnarray}
for all~\(n\) large, where the second inequality in~(\ref{tata}) is obtained by using\\\(\Delta_0 \leq \Delta+1 \leq \beta \cdot n+1\) (see~(\ref{del_not}) and the first condition in~(\ref{m_cond2})). From the statement of Theorem~\ref{thm_comp} and using~\(K = 5,\) we have that
\[\beta + 4K\alpha^2  = \beta + 20 \alpha^2 < \frac{1}{2}\] strictly and so the degree of any vertex in~\(G_z\) is strictly less than~\(\frac{n}{2}\) and also, the sum of degrees of any two vertices in~\(G_z\) is at most~\[(2\beta+40\alpha^2)\cdot n+ 3 + 602 n^{\frac{3}{4}} \cdot \log{n} + \frac{12\sqrt{n}}{10}    < n-10\] for all~\(n\) large. Thus the graph~\(G_z\) can be obtained by the above probabilistic method as in the discussion following~(\ref{del_a}).~\(\qed\)

\renewcommand{\theequation}{A.\arabic{equation}}
\setcounter{equation}{0}
\section*{Appendix}
Throughout we use the following deviation estimate. Let~\(Z_i, 1 \leq i \leq t\) be independent Bernoulli random variables satisfying~\[\mathbb{P}(Z_i = 1) = p_i = 1-\mathbb{P}(Z_i = 0).\] If~\(W_t = \sum_{i=1}^{t} Z_i\) and~\(\mu_t = \mathbb{E}W_t,\) then for any~\(0 < \epsilon < \frac{1}{2}\) we have that
\begin{equation}\label{conc_est_f}
\mathbb{P}\left(\left|W_t-\mu_t\right| \geq \epsilon \mu_t\right) \leq 2\exp\left(-\frac{\epsilon^2}{4}\mu_t\right).
\end{equation}
For a proof of~(\ref{conc_est_f}), we refer to Corollary~\(A.1.14,\) pp.~\(312,\) Alon and Spencer (2008).


%





\subsection*{Acknowledgement}
I thank Professors V. Raman, C. R. Subramanian and the referees for crucial comments that led to an improvement of the paper. I also thank IMSc for my fellowships.


\bibliographystyle{plain}

\end{document}